\documentclass[a4paper,12pt]{article}
\usepackage{geometry,latexsym,amssymb}
\usepackage[latin5]{inputenc}
\usepackage{enumerate}
\usepackage{enumitem}
\usepackage[usenames]{color}
\usepackage{graphicx}
\usepackage{enumerate}
\usepackage{rotating}
\usepackage{multirow}
\usepackage{cite}
\usepackage{amsthm}
\usepackage{tikz}
\usetikzlibrary{matrix,arrows}
\usepackage{amsmath,amscd}
\usepackage{calc}
\usepackage{ifthen}
\usepackage[hang,labelsep=period]{caption}
\usepackage{fancyhdr}

\newtheorem{Def}{Definition}[section]
\newtheorem{Exam}{Example}[section]
\newtheorem{Prop}{Proposition}[section]
\newtheorem{Theo}{Theorem}[section]
\newtheorem{Lem}{Lemma}[section]
\newtheorem{Rem}{Remark}[section]
\newtheorem{Cor}{Corollary}[section]
\newenvironment{Prf}{{\bf Proof:} }{\hfill $\Box$
\mbox{}}

\def\emph{\textbf}

\def\NSCM/(A,B,\mu){\mathsf{NSCM/(A,B,\mu)}}
\def\NSGGd/G{\mathsf{NSGGd/G}}
\def\leq{\leqslant}
\def\geq{\geqslant}

\def\INSGd/G{\operatorname{INSGd/G}}
\def\NSCM/(A,B,\mu){\operatorname{NSCM/(A,B,\mu)}}

\numberwithin{equation}{section} 

\begin{document}

\begin{center}
{\Large \textbf{Generalized Burnside Algebra of type $B_{n}$}} \vspace*{0.5cm}

Hasan Arslan*, Himmet Can

\end{center}

\vspace*{-0.3cm}
*{\small {\textit{Department of  Mathematics, Faculty of Science, Erciyes University, 38039, Kayseri, Turkey, hasanarslan@erciyes.edu.tr}}}\\[0pt]
\begin{center}

\vspace*{-0.1cm}

\textbf{Abstract}
\end{center}

We define the generalized Burnside algebra $HB(W_{n})$ for $B_{n}$-type Coxeter group $W_{n}$ and construct an surjective algebra morphism between Mantaci-Reutenauer algebra ${\sum}'(W_{n})$ and $HB(W_{n})$. Then, by obtaining the primitive idempotents $(e_{\lambda})_{\lambda \in \mathcal{DP}(n)}$ of $HB(W_{n})$, we consider the image $\textrm{res}_{W_{A}}^{W_{n}}e_{B}$ and $\textrm{ind}_{W_{A}}^{W_{n}}e_{B}^{A}$ under restriction and induction map between generalized Burnside algebras. We give an alternative formula to compute the elements number of conjugate classes of $W_{n}$. We also obtain an effective method to determine the size of $\mathcal{C}(S_{n})$ which is the set of elements of type $S_{n}$.

\textbf{Keywords:} Mantaci-Reutenauer Algebra, Generalized Burnside Algebra, Orthogonal Primitive Idempotents

\vspace*{-0.1cm}

\section{\normalsize{INTRODUCTION}}
Let $W_{n}$ be the Coxeter group of type $B_{n}$. In this paper, we first introduce the Burnside algebra generated by isomorphism classes of reflection subgroups of $W_{n}$ corresponding to signed compositions of $n$. It is called \textit{generalized Burnside algebra} of type $B_{n}$. Then we investigate the representation theory of this algebra and its relation to Mantaci-Reutenauer algebra constructed in [11]. In [1], authors have constructed parabolic Burnside algebra by using only standard parabolic subgroups of any finite Coxeter groups. Since the class of the reflection subgroups of $W_{n}$ corresponding to signed compositions of $n$ also contains standard parabolic subgroups, the Burnside algebra we have constructed in Section 3 is more general than parabolic Burnside algebra of type $B_{n}$. We have also constructed the primitive idempotents of generalized Burnside algebra in a parallel manner to [1] and then we have obtained their image under restriction and induction map between generalized Burnside algebras. We follow certain notations in [2] and [3] in denoting some concepts.

\section{\normalsize{NOTATION, PRELIMINARIES}}

\subsection{ Hyperoctahedral group.}

Let $(W_n,S_n)$ denote a Coxeter group of type $B_{n}$ and write its generating set as $S_n=\{t, s_{1},\cdots, s_{n-1}\}$. $W_{n}$ acts by the permutation on the set $X_{n}=\{ -n,\cdots, -1, 1, \cdots, n \}$ such that for every $i \in X_{n}$, $w(-i)=-w(i)$. So we have,
\[
W_{n}=\{w \in \textrm{Perm}(X_{n})~:~ \forall i \in X_{n},~w(-i)=-w(i)\}.
\]

If $J \subset S_{n}$, we will denote by $W_{J}$ \textit{standard parabolic subgroup} of $W_{n}$, which is generated by $J$. A \textit{parabolic subgroup }of $W_{n}$ is a subgroup of $W_{n}$ conjugate to some $W_{J}$. Let $t_{1}:=t$ and $t_{i}:=s_{i-1}t_{i-1}s_{i-1}$ for each i, $2 \leq i \leq n$. Put $T_{n}:=\{t_{1},\cdots,t_{n}\}$. It is well-known that there are the following relations between the elements of $S_{n}$ and $T_{n}$:
\begin{enumerate}
  \item $t_{i}^2=1, s_{j}^2=1~\textit{for all}~ i, j, ~1 \leq i \leq n,~1 \leq j \leq n-1;$
  \item $ts_{1}ts_{1}=s_{1}ts_{1}t;$
  \item $s_{i}s_{i+1}s_{i}=s_{i+1}s_{i}s_{i+1}$~$(1 \leq i \leq n-2);$
  \item $ts_{i}=s_{i}t,~~i>1;$
  \item $s_{i}s_{j}=s_{j}s_{i}$ \textit{for} $|i-j|>1;$
  \item $t_{i}t_{j}=t_{j}t_{i}$ \textit{for} $1\leq i, j\leq n$.
\end{enumerate}
We denote by $l:W_{n}\rightarrow \mathbb{N}$ the length function attached to $S_{n}$. Let $\mathcal{T}_{n}$ denote the reflection subgroup of $W_{n}$ generated by $T_{n}$. It is also clear that $\mathcal{T}_{n}$ is a normal subgroup of $W_{n}$. Now let $S_{-n}=\{s_{1},\cdots, s_{n-1}\}$ and let $W_{-n}$ denote the reflection subgroup of $W_{n}$ generated by $S_{-n}$, where $W_{-n}$ is isomorphic to the symmetric group $\Xi_{n}$ of degree $n$. Thus $W_{n}=W_{-n}\ltimes \mathcal{T}_{n}$. Therefore, we have $|W_{n}|=2^{n}.n!$.

Let $\{e_{1} ,\cdots, e_{n}\}$ be the canonical basis of the Euclidian space ${\mathbb{R}}^{n}$ over $\mathbb{R}$. Let
\begin{equation*}
    \Psi_{n}^+=\{e_{i} : 1 \leq i \leq n \} \cup \{e_{j}+ \lambda e_{i} : \lambda \in \{-1, 1\}~ \textrm{and}~ 1 \leq i < j \leq n \}.
\end{equation*}
Then $\Psi_{n}$ is a root system of $B_{n}$-type. Because of this, we get $t_{j}=s_{e_{j}},~(1\leq j\leq n)$ and $s_{i}=s_{e_{i+1}-e_{i}}, ~(1\leq i\leq n-1)$. The set $\Pi_{n}=\{ e_{1}, e_{2}-e_{1}, \cdots, e_{n}-e_{n-1} \}$ is a simple system of $\Psi_{n}$. The generating set $S_{n}$ of $W_{n}$ is denoted by $\{s_{\alpha} : \alpha\in \Pi_{n}\}$ as a set of simple reflections. To have further information about the Coxeter groups of type $B_{n}$, the readers can see [9] and [10]. We set
$$S_{n}^{'}=S_{n} \cup T_{n}=S_{n}\cup \{t_{1}, \dots, t_{n}\}.$$

A \textit{signed composition} of $n$ is an expression of $n$ as a finite sequence $A=(a_{1}, \cdots, a_{k})$ whose each part consists of non-zero integers such that $\sum_{i=1}^{k}|a_{i}|=n$. We set $|A|=\sum_{i=1}^{k}|a_{i}|$. If exactly $k$ parts appear in $A$, a signed composition of $n$, we say that $A$ has \textit{$k$-length}. The number of $k$ will be denoted by $lg(A)$. We write $\mathcal{SC}(n)$ to denote the set of signed compositions of $n$. We also note that $|\mathcal{SC}(n)|=2.3^{n-1}$. Let $A \in \mathcal{SC}(n)$. We denote by $lg^{+}(A)$ and $lg^{-}(A)$ the length of positive and negative parts of $A$, respectively.

Let $A=(a_{1}, \cdots, a_{k}) \in \mathcal{SC}(n)$. $A$ is said to be \textit{positive} if $a_{i} > 0$ for every $i \geq 1$. We say that $A$ is \textit{parabolic} $a_{i} < 0$ for every $i\geq 2$. If $a_{i} \geq -1$ for every $i \geq 1$ then we called $A$ as \textit{semi-positive} . Let define $A^{+}=(|a_{1}|,\cdots, |a_{r}|)$. Then $A^{+}$ is a positive signed composition of $n$. The set of positive and parabolic signed compositions of $n$ is denoted by $\mathcal{SC}^{+}(n)$ and $\mathcal{SC}_{p}(n)$, respectively. If $A=(a_{1}, \cdots, a_{r})\in \mathcal{SC}(m)$ and $B=(b_{1}, \cdots, b_{s})\in \mathcal{SC}(n)$, then $A \sqcup B$ will denote the signed composition $(a_{1}, \cdots, a_{r},b_{1}, \cdots, b_{s})$ of $m+n$.

A \textit{double partition} $\mu=({\mu}^{+}, {\mu}^{-})$ of $n$ consists of a pair of partitions ${\mu}^{+}$ and ${\mu}^{-}$ such that $|\mu|=|{\mu}^{+}|+|{\mu}^{-}|=n$. We denote the set of double partitions of $n$ by $\mathcal{DP}(n)$. For $\mu=({\mu}^{+}, {\mu}^{-}) \in \mathcal{DP}(n)$ we set $\hat{\mu}:=\mu^{+} \sqcup \mu^{-}$, then $\hat{\mu}\in \mathcal{SC}(n)$. Because of this description, the map $\mathcal{DP}(n) \rightarrow \mathcal{SC}(n)$,~$\mu \mapsto \hat{\mu}$ is injective.

Now let $A \in \mathcal{SC}(n)$. If $\mu^{+}$ (resp. $\mu^{-}$) is rearrangement of the positive parts (resp. negative parts)of $A$ in increasing order, then $\boldsymbol{\lambda}(A):=(\lambda^{+}, \lambda^{-})$ is a double partition of $n$. Thus, the map $\boldsymbol{\lambda} : \mathcal{SC}(n) \rightarrow \mathcal{DP}(n)$ is surjective, and furthermore $ \boldsymbol{\lambda}(\hat{\mu})=\mu$.

In [2], Bonnafé ve Hohlweg have constructed some reflection subgroups of $W_{n}$ corresponding to signed compositions of $n$ as an analogue of $\Xi_{n}$. For each $A=(a_{1}, \cdots, a_{k})\in \mathcal{SC}(n)$, the reflection subgroup $W_{A}$ of $W_{n}$ is generated by $S_{A}$, which is
\begin{align*}
    S_{A}= & \{ s_{p}\in W_n~:~|a_{1}|+ \cdots + |a_{i-1}|+1 \leq p \leq |a_{1}|+ \cdots + |a_{i}|-1\} \\
           & \cup \{ {t_{|a_{1}|+ \cdots + |a_{j-1}|+1}} \in T_{n} \} ~|~ a_{j} > 0 \} \subset S_{n}^{'}
\end{align*}
By the definition of $S_{A}$, we have the isomorphism $W_{A}\cong W_{a_{1}} \times \cdots \times W_{a_{k}}$. By taking into account the definition of the generating set $S_{A}$ and the isomorphism $W_{A}\cong W_{a_{1}} \times \cdots \times W_{a_{r}}$, if $a_{i}>0$ for $i,~1 \leq i \leq r$ then $\textrm{rank} W_{a_{i}}=a_{i}$ or if $a_{i}<0$ then $\textrm{rank} W_{a_{i}}=|a_{i}|-1$. Therefore, we have
\begin{equation*}
    \textrm{rank}W_{A}=|S_{A}|=n-lg^{-}(A).
\end{equation*}
Because of $\sum_{i=1}^{r}|a_{i}|=n$, we obtain $\textrm{rank}W_{A}=|S_{A}| \leq n$.

Let $S_{A}^{'}=S_{n}^{'} \cap W_{A},~~\Psi_{A}=\{\alpha \in \Psi_{n} : s_{\alpha}\in W_{A}\}~\textrm{and}~ \Psi_{A}^+=\Psi_{A} \cap \Psi_{n}^+$. In addition, $\Psi_{A}^+$ is a positive root system of $\Psi_{A}$ and $\Pi_{A}$ is a simple system of $\Psi_{A}$ contained in $\Psi_{A}^+$. Thus $S_{A}=\{s_{\alpha} : \alpha \in \Pi_{A}\}$ and so by [2] $(W_{A}, S_{A})$ is a Coxeter system. If $T_{A}$ is defined as $T_{n}\cap W_{A}$, then $\mathcal{T}_{A}=\mathcal{T}_{n} \cap W_{A}=\langle T_{A} \rangle$. Thus, there exists a semidirect decomposition as follows:

$$W_{A}=\langle W_{A} \cap S_{-n} \rangle \ltimes \mathcal{T}_{A}.$$

For $A,B \in \mathcal{SC}(n)$, we write $A\subset B$ if $W_{A} \subset W_{B}$, where $\subset$ is a partial ordering relation on $\mathcal{SC}(n)$. Note that $S_{A}^{'}=S_{B}^{'}$ if and only if $A=B$. Let $A\in \mathcal{SC}(n)$. Denote by $\textit{c}_{A}$ a Coxeter element of $W_{A}$ in terms of generating set $S_{A}$. For $B, B^{'} \subset A$, we write $B \equiv_{A} B{'}$ if $W_{B}$ is conjugate to $W_{B^{'}}$ under $W_{A}$. Also, a necessary and sufficient condition that $\textit{c}_{B}$ and $\textit{c}_{B^{'}}$ are conjugate to each other in $W_{A}$ is that $B \equiv_{A} B^{'}$. This equivalence is a special case for these kind of reflection subgroups of $W_{n}$, because this statement is not true for every reflection subgroups of $W_{n}$. By [7], although some two reflection subgroup $R$ and $R'$ of $W_{n}$ contain $W_{n}$-conjugate Coxeter element $c$ and $c'$ respectively, they can not be $W_{n}$-conjugate to each other. Every element of $W_{n}$ is $W_{n}$-conjugate to $\textit{c}_{A}$ for some $A\in \mathcal{SC}(n)$. Bonnafé [3] has showed that for $A, B \in \mathcal{SC}(n)$, $W_{A}$ is conjugate to $W_{B}$ in $W_{n}$ if and only if $\boldsymbol{\lambda}(A)=\boldsymbol{\lambda}(B)$. The number of conjugate classes of $W_{n}$ is equal to $|\mathcal{DP}(n)|$ and so we may split up $W_{n}$ into $|\mathcal{DP}(n)|$ equivalence classes.

\subsection{ Mantaci-Reutenauer algebra.}

For any $A\in \mathcal{SC}(n)$, we set
\begin{equation*}
              D_{A}=\{x\in W_{n}~:~\forall ~s\in S_{A},~~l(xs)>l(x) \}.
\end{equation*}
By [2] and [9], $D_{A}$ is the distinguished coset representatives of $W_{A}$ in $W_{n}$ provided that the map $D_{A}\times W_{A}\rightarrow W_{n},~(d,w)\mapsto dw $ is bijective. For $A, B\in \mathcal{SC}(n)$ such that $B\subset A$, the set $D_{B}^A=D_{B}\cap W_{A}$ is the distinguished coset representatives of $W_{B}$ in $W_{A}$.
Let
$$d_{A}=\sum_{w\in D_{A}}w \in \mathbb{Q}W_{n},$$
and let
$${\sum}'(W_{n})=\bigoplus_{A \in \mathcal{SC}(n)}\mathbb{Q}d_{A}.$$
For every $A\in \mathcal{SC}(n)$, $\theta_{n}~:~{\sum}'(W_{n}) \rightarrow \mathbb{Q}\textrm{Irr}W_{n}$ be the unique $\mathbb{Q}$-linear map such that $\theta_{n}(d_{A})=\textrm{Ind}_{W_{A}}^{W_{n}}1_{A}$, where $\mathbb{Q}\textrm{Irr}W_{n}$ and $1_{A}$ stands for the algebra of irreducible characters of $W_{n}$ and the trivial character of $W_{A}$, respectively. For $A, B \in \mathcal{SC}(n)$, we put
$$D_{AB}=D_{A}^{-1} \cap D_{B}.$$
For a $x\in D_{AB}$, by [2],
$$W_{A}\cap {{^{x}}W_{B}}=W_{A\cap {{^{x}}B}}$$
and $x$ is the unique element of $W_{A}xW_{B}$ of minimal length.

\begin{Theo}[2]
If $A, B$ are two signed composition of $n$, then
\begin{enumerate}
  \item [(a)] ${\sum}'(W_{n})$ is a $\mathbb{Q}$-subalgebra of $\mathbb{Q}W_{n}$;
  \item [(b)] $\theta_{n}$ is a surjective morphism of $\mathbb{Q}$-algebras;
  \item [(c)] $Ker\theta_{n}=\bigoplus_{{A, B \in \mathcal{SC}(n)},~{A\equiv_{n}B'}}\mathbb{Q}(d_{A}-d_{B})$;
  \item [(d)] $Ker\theta_{n}$ is the radical of ${\sum}'(W_{n})$.
\end{enumerate}
\end{Theo}
${\sum}'(W_{n})$ is called \textit{Mantaci-Reutenauer algebra} of $W_{n}$, which includes the classical Solomon algebras of $W_n$ and $\Xi_{n}$. For a subset $\mathcal{F}$ of $\mathcal{SC}(n)$, we set
\[
{\sum}_{\mathcal{F}}'(W_{n})= \bigoplus_{A \in \mathcal{F}}\mathbb{Q}d_{A}.
\]
Bonnafé~\cite{Br3} has introduced the order $\preceq $ on $\mathcal{SC}(n)$, which is more useful than $\subset$ to determine the structure of the multiplication in ${\sum}'(W_{n})$ in the following way. For $A, B \in \mathcal{SC}(n)$, we write $A \preceq B$ if and only if either
\begin{itemize}
  \item $A \subset B$ or
  \item $A \subset B^+$,~$lg(A) > lg(B)$,~ $lg^{-}(A)\geq lg^{-}(B)$
\end{itemize}
satisfies.
Let $\mathcal{F}_{\preceq B}=\{A\in \mathcal{SC}(n)~:~ A\preceq B \}$ and let ${\sum}_{\preceq B}'(W_{n})={\sum}_{\mathcal{F}_{\preceq B}}'(W_{n})$.
For $A$ and $B$ be two signed composition of $n$, we define
$$A\subset_{\boldsymbol{\lambda}} B ~\Leftrightarrow~\boldsymbol{\lambda}(A) \subset \boldsymbol{\lambda}(B).$$
Here the relation $\subset_{\boldsymbol{\lambda}}$ is an reflexive and transitive on $\mathcal{SC}(n)$. In other words, $A\subset_{\boldsymbol{\lambda}} B$ if and only if $W_{A}$ is contained in some conjugate of $W_{B}$. We will write $A\subsetneq_{\boldsymbol{\lambda}} B$ if $\boldsymbol{\lambda}(A) \subsetneq \boldsymbol{\lambda}(B)$. We set $\mathcal{F}_{{\subset}_{\boldsymbol{\lambda}} B}=\{A\in \mathcal{SC}(n)~:~ A{\subset}_{\boldsymbol{\lambda}} B \}$ and as an abbreviation we shall simply ${\sum}_{{\subset}_{\boldsymbol{\lambda}} B}'(W_{n})$ for ${\sum}_{\mathcal{F}_{{\subset}_{\boldsymbol{\lambda}} B}}'(W_{n})$. We will use the following lemma proved by Bonnafé [3] in the beginning of construction of generalized Burnside algebra.

\begin{Prop}[3]
Let $A$ and $B$ be any two signed composition of $n$. Then,
\begin{enumerate}
\item[\textbf{\textrm{(a)}}] There is a map $f_{AB} : D_{AB} \rightarrow \mathcal{SC}(n)$ satisfying the following conditions:
\begin{itemize}
\item For every $x \in D_{AB}$, $f_{AB}(x)\subset B$ and $f_{AB}(x)\equiv_{B} {^{x^{-1}}A} \cap B$.
\item $d_{A}d_{B}-\sum _{x\in D_{AB}}d_{f_{AB}(x)} \in {\sum}_{{\subsetneq}_{\boldsymbol{\lambda}} A}'(W_{n}) \cap {\sum}_{\prec B}'(W_{n}) \cap Ker \theta_{n}$.
\end{itemize}
\item[\textbf{\textrm{(b)}}] If A parabolic or $B$ is semi-positive, then $f_{AB}(x)={^{x^{-1}}A} \cap B$ for $x \in D_{AB}$ and
$d_{A}d_{B}=\sum _{x\in D_{AB}}d_{{^{x^{-1}}A} \cap B}$.
\end{enumerate}
\end{Prop}

For $\lambda\in \mathcal{DP}(n)$ we set $\tau_{\lambda}~:~ {\sum}'(W_{n})\rightarrow \mathbb{Q},~x\mapsto \theta_{n}(x)(c_{\lambda})$. The map $\tau_{\lambda}$ is independent of the choice of $c_{\lambda}\in C(\lambda) $ and it is also an algebra morphism. Moreover, the set $\{\tau_{\lambda}~:~\lambda\in \mathcal{DP}(n)\}$ is the collection of irreducible representations of ${\sum}'(W_{n})$.

Throughout this paper, we will use the following facts proved in [3] frequently .

\begin{Prop}[3]
For $A, B\in \mathcal{SC}(n)$, let define the sets $D_{AB}^{\subset}=\{x\in D_{AB} :~^{x^{-1}}W_{A}\subset W_{B}\}$ and $D_{AB}^{\equiv}=\{x \in D_{AB} : W_{A}=^{x}W_{B}\}$. Then the following statements hold:
\begin{enumerate}
\item $\tau_{\boldsymbol{\lambda}(A)}(d_{B})=|D_{AB}^{\subset}|$.
\item $D_{AB}^{\equiv}=\{x\in W_{n}~:~\Pi_{A}=^{x}\Pi_{B}\}$.
  \item $\mathcal{W}(B)=\{w\in W_{n} : w(\Pi_{B})=\Pi_{B}\}$;
  \item $\mathcal{W}(B)$ is a subgroup of $N_{{W_{n}}}(W_{B})$;
  \item $N_{{W_{n}}}(W_{B})=\mathcal{W}(B)\ltimes W_{B}$.
\end{enumerate}
\end{Prop}
If $A\equiv B$, then it is clear $D_{AB}^{\equiv}=D_{AB}^{\subset}$. For $A\in \mathcal{SC}(n)$ we set $\mathcal{W}(A)=D_{AA}^{\subset}$.

\section{Generalized Burnside Algebra of \textbf{$W_{n}$}}

Let $A, B$ be any two signed composition of $n$. By [2], we have that
$$A\equiv_{n}B~ \Leftrightarrow ~W_{A}\sim W_{B}~\Leftrightarrow ~[W/W_{A}]=[W/W_{B}]$$
where $[W/W_{A}]$ represents the isomorphism class of $W_{n}$-set $W/W_{A}$. The orbits of $W_{n}$ on $W/W_{A}\times W/W_{B}$ are of the form $(W_{A}x,W_{B})$ where $x\in D_{AB}$. The stabilizer of $(W_{A}x,W_{B})$ in $W_{n}$ is  $^{x^{-1}}W_{A}\cap W_{B}=W_{^{x^{-1}}A\cap B}$. Therefore
$$[W/W_{A}].[W/W_{B}]=[W/W_{A} \times W/W_{B}]=\sum_{x\in D_{AB}}[W/W_{^{x^{-1}}A\cap B}].$$
Thus, we are now in a position to give the following definition.

\begin{Def} \normalfont
\textit{The generalized Burnside algebra} of $W_{n}$ is $\mathbb{Q}$-spanned by the set $\{ [W/W_{A}] : A \in \mathcal{SC}(n) \}$ and it is denoted by $\textrm{HB}(W_{n})$.
\end{Def}

From Proposition 2.1 and the structure of $\textrm{Ker}(\theta_{n})$, the multiplication in ${\sum}'(W_{n})$ may be expressed by
$$d_{A}d_{B}=\sum _{x\in D_{AB}}d_{f_{AB}(x)} + {\sum}_{N\equiv_{n}N'}a_{NN'}(d_{N}-d_{N'}),$$
where $a_{NN'}\in \textrm{Z}$; $N, N'{\subsetneq}_{\boldsymbol{\lambda}} A$; $N, N'\prec B$; $f_{AB}(x)\subset B$ and $f_{AB}(x)\equiv_{B} {^{x^{-1}}A} \cap B$.

Now we define
$$\psi~:~{\sum}'(W_{n})\rightarrow HB(W_{n}),~d_{A}\mapsto [W/W_{A}].$$
the map $\psi$ is well-defined and surjective. On considering the structure of $\textrm{Ker}\theta_{n}$ and $f_{AB}(x)\equiv_{B} {^{x^{-1}}A} \cap B \Rightarrow W_{f_{AB}(x)} \sim_{W_{B}} W_{^{x^{-1}}A \cap B}$, we get
\begin{align*}
\psi (d_{A}d_{B})& =\psi(\sum _{x \in D_{AB}} d_{f_{AB}(x)} + {\sum}_{N \equiv_{n} N'} a_{NN'} (d_{N}-d_{N'})) \\
& =\sum_{x\in D_{AB}} \psi(d_{f_{AB}(x)}) + {\sum}_{N\equiv_{n}N'}a_{NN'} (\psi(d_{N})-\psi (d_{N'}))  \\
& =\sum_{x\in D_{AB}} [W/W_{f_{AB}(x)}] \\
& =\sum_{x\in D_{AB}} [W/{W_{^{x^{-1}}A \cap B}}] \\
& =[W/W_{A}].[W/W_{B}] \\
& =\psi (d_{A})\psi (d_{B}).
\end{align*}
Then the map $\psi$ is an algebra morphism. Since $\textrm{dim}_{\mathbb{Q}}\textrm{HB}(W_{n})=\textrm{dim}_{\mathbb{Q}} \mathbb{Q}\textrm{Irr}W_{n}=|\mathcal{DP}(n)|$, then there is an algebra isomorphism between  $\mathbb{Q}\textrm{HB}(W_{n})$ and $\mathbb{Q}\textrm{Irr}W_{n}$ such that

$$ \textrm{HB}(W_{n}) \rightarrow \mathbb{Q}\textrm{Irr}W_{n},~[W/W_{A}]\mapsto \textrm{ind}_{W_{A}}^{W_{n}}1_{A}.$$

Now let $\lambda, \mu \in \mathcal{DP}(n)$ and let $\varphi_{\lambda}=\textrm{ind}_{W_{A}}^{W_{n}}1_{A}$ for any $A\in \boldsymbol{\lambda}^{-1}(\lambda)$. By [2], $\varphi_{\lambda}(c_{\lambda})=\pi_{\lambda}(x_{\hat{\lambda}})=|D_{\hat{\lambda}\hat{\lambda}}^{\subset}|\neq 0$ and $\pi_{\lambda}(x_{\hat{\mu}})=0$ if $\lambda \nsubseteq \mu$. Thus the matrices $(\pi_{\lambda}(c_{\hat{\mu}}))_{\lambda, \mu \in \mathcal{DP}(n)}$ is lower diagonal. Then $(\varphi_{\lambda}(c_{\mu}))_{\lambda, \mu}$ is upper diagonal and has also positive diagonal entries. Therefore $(\varphi_{\lambda}(c_{\mu}))_{\lambda, \mu}$ is invertible and its inverse is $(u_{\lambda \mu})_{\lambda, \mu \in \mathcal{DP}(n)}$. As in [1], we define
\begin{equation}
    e_{\lambda}=\sum_{\mu \in \mathcal{DP}(n)}u_{\lambda \mu}\varphi_{\mu}.
\end{equation}
By definition of $e_{\lambda}$ and $(\varphi_{\lambda}(c_{\mu}))^{-1}=(u_{\lambda \mu})$, we obtain that
\[
e_{\lambda}(c_{\mu})=\sum_{\gamma \in \mathcal{DP}(n)}u_{\lambda \gamma}\varphi_{\gamma}(c_{\mu})=\delta_{\lambda, \mu}.
\]
Hence the set $\{e_{\lambda} : \lambda \in \mathcal{DP}(n)\}$ is a collection of primitive idempotents of $\textrm{HB}(W_{n})$. Wedderburn's structure theorem which we need in our works is given as follows:

\begin{Theo}[5]
Let $A$ be an algebra over the field $K$ and let $J(A)$ be the Jacobson radical of A. Then we have
\begin{equation*}
    A/J(A)\cong {\oplus_{i=1}^{s}} A_{i},
\end{equation*}
where $A_{i}$ is an matrix algebra isomorphic to $Mat_{n_{i}}(K)$.
\end{Theo}

\begin{Lem}[5]
Let $A$ be an commutative algebra over the field $K$. Every algebra map $A \rightarrow K$ is of the form $s_{i}$ such that $s_{i}(1_{A_{j}})=\delta_{i,j}$, where the $1_{A_{i}}$'s are primitive idempotents.
\end{Lem}

For each $A \in \mathcal{SC}(n)$ by \cite{Br6},
$$s_{A} : HB(W_{n})\rightarrow \mathbb{Q},~s_{A}([X])=|^{W_{A}}X|$$
is an algebra map, where $^{W_{A}}X=\{x\in X : W_{A}x=x\}$. Since $HB(W_{n})$ is semisimple and commutative algebra, then all the maps $HB(W_{n})\rightarrow \mathbb{Q}$ are of the form $s_{A}$ for every $A\in \mathcal{SC}(n)$.
Consequently, we have $HB(W_{n})=\oplus_{\lambda\in \mathcal{DP}(n)}\mathbb{Q}e_{\lambda}$. The proof of the following lemma is immediately seen from [6].

\begin{Lem}
For $A, B \in \mathcal{SC}(n)$, we have that
$$s_{A}=s_{B}~\Leftrightarrow ~ \boldsymbol{\lambda}(A)=\boldsymbol{\lambda}(B).$$
\end{Lem}

The dual basis of $HB(W_{n})$ is $\{s_{\boldsymbol{\lambda}(A)} : \lambda \in \mathcal{DP}(n) \}$. For any $\lambda,~\mu \in \mathcal{DP}(n)$, there exists the following equality
\begin{equation}
    s_{\lambda}(e_{\mu})=\delta_{\lambda, \mu},
\end{equation}
and so any element $x$ in $HB(W_{n})$ can be expressed as $x=\sum_{\lambda \in \mathcal{DP}(n)}s_{\lambda}(x)e_{\lambda}$.

Let $A$ be a signed composition of $n$. Induction and restriction of characters give rise to two maps between $HB(W_{A})$ and $HB(W_{n})$. For any $A, B\in \mathcal{SC}(n)$ such that $B \subset A$, we have $\textrm{ind}_{W_{A}}^{W_{n}}([W_{A}/W_{B}])=[W_{n}/W_{B}]$.

\begin{Def} \normalfont
Let be $A, B\in \mathcal{SC}(n)$ such that $B \subset A$. The \textit{restriction} is a linear map
$$ \textrm{res}_{W_{B}}^{W_{A}}:  HB(W_{A}) \rightarrow  HB(W_{B}),~ \textrm{res}_{W_{B}}^{W_{A}}([W_{A}/W_{C}])=\sum_{d \in W_{A}\cap D_{CB}}[W_{B}/W_{B\cap ^{d^{-1}}C}].$$
\end{Def}

Before going into a further discussion of the restriction and induced character theories of generalize Burnside algebra, we shall first mention the number of elements of the conjugacy class of $W_{A}$ in $W_{n}$.

Let $A, B \in \mathcal{SC}(n)$. We define
$$\textrm{inv}_{(W_{n}/W_{A})}(W_{B})=\{xW_{A} : x \in W_{n},~W_{B}xW_{A}=xW_{A}\},$$
which is the set of elements in $W_{n}/W_{A}$ fixed by $W_{B}$. Also $\textrm{inv}_{W_{n}/W_{A}}(W_{B})=\{xW_{A} : x \in W_{n},~x^{-1}W_{B} x\leq W_{A}\}$ and so $\textrm{inv}_{W_{n}/W_{A}}(W_{B})= \varnothing $ unless $\boldsymbol{\lambda}(B) \not \subset \boldsymbol{\lambda}(A)$. The \textit{mark} of $W_{B}$ on $W_{n}/W_{A}$ is defined to be the number
$$\mu_{AB}=|\textrm{inv}_{W_{n}/W_{A}}(W_{B})|.$$
From [6], for each $A'\in \boldsymbol{\lambda}(A)$ and $B' \in \boldsymbol{\lambda}(B)$, it is clear $\mu_{A'B'}=\mu_{AB}$.

\begin{Def} \normalfont
The \textit{signed table of marks} of $W_{n}$ is a square matrix of type $|\mathcal{DP}(n)|\times|\mathcal{DP}(n)|$ as follows
$$M(W_n)=(\mu_{\lambda \vartheta})_{\lambda, \vartheta\in \mathcal{DP}(n)}.$$
\end{Def}

This table also contains the parabolic table of marks of $W_{n}$. $M(W_n)$ is lower triangular with non-zero diagonal entries and so it is invertible. Since $|D_{BA}^{\subset}|$ is the number of fixed points in $W_{n}/W_{A}$ under the action of $\textrm{c}_B$, we have
$$\mu_{AB}=|\textrm{inv}_{W_{n}/W_{A}}(W_{B})|=|D_{BA}^{\subset}|=\pi_{\boldsymbol{\lambda}(B)}(d_{A}).$$
As a corollary, each entries of the matrix $M(W_n)$ can also be read from the matrices $(\pi_{\lambda}(c_{\hat{\mu}}))_{\lambda, \mu \in \mathcal{DP}(n)}$.

\begin{Prop}
Let $A, B \in \mathcal{SC}(n)$. Then we have
$$\mu_{AB}=|\textrm{inv}_{W_{n}/W_{A}}(W_{B})|=|\{^{x}W_{A} : x \in W_{n}, W_{B}\leq ^{x}W_{A}\}| \cdot |\mathcal{W}(A)|.$$
\end{Prop}
The proof is an easy application of the definition of $\mu_{AB}$.
As a result of the above proposition, the number of conjugates of $W_{A}$ containing $W_{B}$ is
$$|\{^{x}W_A : x\in W_{n}, W_{B}\leq ^{x}W_A \}|=\frac{|D_{BA}^{\subset}|}{|\mathcal{W}(A)|}.$$

\begin{Prop}
Let $A \in \mathcal{SC}(n)$ and $\boldsymbol{\lambda}(A)= \lambda $. The number of all reflection subgroups of $W_{n}$ that are conjugate to $W_{A}$ is
$$|[W_A]|=|D_A| \cdot u_{\lambda, \lambda}.$$
\end{Prop}

\begin{Prf}
Put $[W_A]=\{^{x}W_{A} : x \in W_{n}\}$. Now we note that $xW_{A}x^{-1}=yW_{A}y^{-1}$ if and only if $x \in yN_{W_{n}}(W_{A})$. Thus, the number of distinct conjugates of $W_{A}$ in $W_{n}$ is $[W_{n} : N_{W_{n}}(W_{A})]$. Since also $N_{W_{n}}(W_{A})=\mathcal{W}(A)\ltimes W_{A}$, we have
$$|[W_A]|=\frac{|W_{n}|}{|\mathcal{W}(A)|\cdot|W_{A}|}=\frac{|D_{A}|}{|\mathcal{W}(A)|}.$$
Likewise, from the fact that $\pi_{\boldsymbol{\lambda}(A)}(d_{A})=|D_{AA}^{\subset}|=|\mathcal{W}(A)|$ and $\varphi_{\lambda}(cox_{\lambda})=\pi_{\boldsymbol{\lambda}(A)}(d_{A})=\frac{1}{u_{\lambda,\lambda}}$, we see that $|[W_A]|$ has the desired number.
\end{Prf}

\begin{Exam} \normalfont
The set $D_{(2,1)}=\{1, s_2, s_1s_2\}$ is the distinguished coset representatives of reflection subgroup $W_{(2,1)}$ in $W_3$, then the number of all reflection subgroups conjugate to $W_{(2,1)}$ in $W_{3}$ is
$$|[W_{(2,1)}]|=|D_{(2,1)}| \cdot u_{(2,1),(2,1)}=3 \cdot 1=3.$$
These are explicitly $W_{(2,1)}$, $W_{(1,2)}$ and $^{s_2}W_{(2,1)}=\langle s_2s_1s_2, t_1, t_2 \rangle$. We note that the last one is not a reflection subgroup of $W_{3}$ corresponding to any signed composition of $3$.
\end{Exam}

\begin{Rem}
For $A, B\in \mathcal{SC}(n)$ such that $B \subset A$ and for any $x\in HB(W_{n})$, by using the definition of $s_{A}$ one can see that there is the relation  $s_{B}^{A}(\textrm{res}_{W_{A}}^{W_{n}}(x))=s_{B}(x).$
\end{Rem}

We can now five the following proposition which is analog to Theorem 3.2.4 in [1].

\begin{Prop}
Let be $A, B\in \mathcal{SC}(n)$ and let $A_{1}, A_{2}, \cdots, A_{r}$ be representatives of the $W_{A}$-equivalent classes of subsets of $A$, which is $W_{n}$-equivalent to $B$. Then,
$$ \textrm{res}_{W_{A}}^{W_{n}}e_{B}={\sum_{i=1}^{r}}e_{B_{i}}^{A}.$$
If $B$ is not $W_{n}$-equivalent to any subset of $A$ then $\textrm{res}_{W_{A}}^{W_{n}}e_{B}=0$.
\end{Prop}

\begin{Prf}
Since $\textrm{res}_{W_{A}}^{W_{n}} e_{B}$ is an element of $\mathbb{Q}HB(W_{A})$ we have $res_{W_{A}}^{W_{n}}e_{B} =\sum_{A_{i}\subset A}s_{B_{i}}^{A}(\textrm{res}_{W_{A}}^{W_{n}}(e_{B}))\xi_{C_{i}}^{C}.$
Then by using Remark (3.1) and the relation (3.2), we get
\begin{align*}
\textrm{res}_{W_{A}}^{W_{n}}e_{B} &=\sum_{A_{i}\subset A}s_{A_{i}}(e_{B})e_{A_{i}}^{A} \\
&=\sum_{\substack{A_{i}\subset A \\ A_{i} \equiv_{A} B}} e_{A_{i}}^{A} \\
&={\sum_{i=1}^{r}}e_{B_{i}}^{A}.
\end{align*}
\end{Prf}

\begin{Prop}
Let $A, B\in \mathcal{SC}(n)$ and let $B\subset A$.Then we have
$$ \textrm{ind}_{W_{A}}^{W_{n}}e_{B}^{A}=\frac{|\mathcal{W}(B)|}{|W_{A}\cap \mathcal{W}(B)|}\cdot e_{B}.$$
\end{Prop}

\begin{Prf}
Firstly, we assume that $A=B$ and $\textit{c}_{A}$ is a Coxeter element of $W_{A}$.
Since the image of $c_{A}$ under permutation character of $W_{n}$ on the cosets of $W_{A}$ is $|\mathcal{W}(A)|$ then it follows from the fact that
$$x^{-1}c_{A}x \in W_{A} \Leftrightarrow x\in N_{W_{n}}(W_{A}).$$
Thus we obtain
\begin{align*}
\textrm{ind}_{W_{A}}^{W_{n}}e_{A}^{A}(c_{A})&=|D_{A}\cap N_{W_{n}}(W_{A})| \\
&=|\mathcal{W}(A)|.
\end{align*}
As $\textrm{ind}_{W_{A}}^{W_{n}}e_{A}^{A}$ takes value zero except for the elements conjugate to $c_{A}$ and so we get
$$\textrm{ind}_{W_{A}}^{W_{n}}e_{A}^{A}=|\mathcal{W}(A)|e_{A}.$$
By transitivity of induced characters, we generally get
\begin{align*}
\textrm{ind}_{W_{A}}^{W_{n}}e_{B}^{A}&=\textrm{ind}_{W_{A}}^{W_{n}}(\frac{1}{|W_{A}\cap \mathcal{W}(B)|}|W_{A}\cap \mathcal{W}(B)|e_{B}^{A}) \\
&= \textrm{ind}_{W_{A}}^{W_{n}}(\frac{1}{|W_{A}\cap \mathcal{W}(B)|}\textrm{ind}_{W_{B}}^{W_{A}}e_{B}^{B}) \\
&=\frac{|\mathcal{W}(B)|}{|W_{A}\cap \mathcal{W}(B)|}e_{B}.
\end{align*}
\end{Prf}

Furthermore, there is also the equality $\textrm{ind}_{W_{A}}^{W_{n}}e_{B}^{A}=|N_{W_{n}}(W_{B}):N_{W_{A}}(W_{B})|e_{B}$.

\begin{Theo}
Let $A, B\in \mathcal{SC}(n)$ and let $\boldsymbol{\lambda}(B) \subset \boldsymbol{\lambda}(A)$. If $B_{1}, B_{2}, \cdots, B_{r}$ are the representatives of the $W_{A}$-equivalent classes of subsets of $A$, $W_{n}$-equivalent to $B$, then for $\textit{c}_{B}\in W_{n}$,
$$ \textrm{ind}_{W_{A}}^{W_{n}}1_{A}(\textrm{c}_{B})=\sum_{i=1}^{r}\frac{|\mathcal{W}(B)|}{|W_{A}\cap \mathcal{W}(B_{i})|}.$$
\end{Theo}

\begin{Prf}
Let $A, B\in \mathcal{SC}(n)$. If $A\equiv_{n}B$ then it is easy to prove that $|\mathcal{W}(A)|=|\mathcal{W}(B)|$. $1_{A}=\sum_{E}e_{E}^{A}$, where $E \in \mathcal{SC}(n)$ runs over $W_{A}$-conjugate classes of subsets of $A$. From Proposition 9, we have
$$\textrm{ind}_{W_{A}}^{W_{n}}1_{A}=\sum_{E}\textrm{ind}_{W_{A}}^{W_{n}}e_{E}^{A} \Rightarrow \textrm{ind}_{W_{A}}^{W_{n}}1_{A}=\sum_{E}\frac{|\mathcal{W}(E)|}{|W_{A}\cap \mathcal{W}(E)|}\cdot e_{E}.$$
As each $B_{i}$ is $W_{n}$-congruent to $B$, then $e_{E}(\textrm{c}_{B})=1$  if and only if $E\equiv_{W_{A}} B_{i}$. Thus from (3.1), we obtain that
\[
\textrm{ind}_{W_{A}}^{W_{n}}1_{A}(\textrm{c}_{B})=\sum_{i=1}^{r}\frac{|\mathcal{W}(B)|}{|W_{A}\cap \mathcal{W}(B_{i})|}.
\]
Hence the theorem is proved.
\end{Prf}

Theorem 3.3 and Proposition 3.5 give us a useful computation to determine the coefficient of the sign character $\varepsilon_{n}$ in the expression of primitive idempotent $e_{\lambda},~\lambda \in \mathcal{DP}(n)$ in terms of irreducible characters of $W_{n}$. From now on, for convention, we simply write $e_{\hat{\lambda}}$ instead of $e_{\lambda}$ to avoid the complexity of parenthesis in some calculation.

\begin{Theo}
$u_{(n),(-1,\cdots, -1)}=\frac{(-1)^n}{2n}$.
\end{Theo}

\begin{Prf}
Let $\varepsilon_{n} : W_{n} \rightarrow \{-1,~ 1\}$ and $\chi_{reg} : W_{n} \rightarrow \mathbb{Z}$ be the sign and regular character of $W_{n}$, respectively. For $A=(-1,\cdots, -1)$ it is satisfied $\textrm{ind}_{W_{A}}^{W_{n}}1_{A}=\chi_{reg}$. The character $\varepsilon_{n}$ is contained in $\chi_{reg}$ with the property that its coefficient is just 1, thus we have
$$\langle \textrm{ind}_{W_{A}}^{W_{n}}1_{A},\varepsilon_{n}\rangle=1.$$
Now let $A \neq (-1,\cdots, -1)$. By Frobenius Reciprocity and $\sum_{w \in W_{A}}(-1)^{l{(w)}}=0$, it is obtained that $\langle \textrm{ind}_{W_{A}}^{W_{n}}1_{A},\varepsilon_{n}\rangle=0.$
If $w$ is conjugate to $c_{W_{n}}$, then $e_{(n)}(w)=1$ and $l(w)=n$. Let $ccl_{W_{n}}(c_{W_{n}})$ denotes the conjugate class of $c_{W_{n}}$ in $W_{n}$. By [4], by considering the fact $|ccl_{W_{n}}(c_{W_{n}})|=\frac{|W_{n}|.n}{2N}$, we have
\[
\langle e_{(n)},\varepsilon_{n}\rangle =\frac{(-1)^{n}}{2n}.
\]
On the other hand, $\langle e_{(n)},\varepsilon_{n}\rangle=\sum_{\mu \in \mathcal{DP}(n)} u_{(n) \mu} \langle \varphi_{\mu}, \varepsilon_{n} \rangle=u_{(n),(-1,\cdots, -1)}$ and so the proof is completed.
\end{Prf}

\begin{Prop}
For $\lambda \in \mathcal{DP}(n)$ and $\lambda \neq (n)$, then we have
$$u_{\lambda,(-1,\cdots, -1)}=(-1)^{|S_{\hat{\lambda}}|}\cdot \frac{|\mathcal{C}(\lambda)|}{|W_{n}|}.$$
\end{Prop}

\begin{Prf}
Since the sign character is constant on conjugate classes, then we have
\begin{align*}
\langle e_{\lambda},\varepsilon_{n}\rangle &=\frac{1}{|W_{n}|}\sum_{w \in \mathcal{C}(\lambda)} (-1)^{l(w)}  \\
&=(-1)^{|S_{\hat{\lambda}}|}\cdot \frac{|\mathcal{C}(\lambda)|}{|W_{n}|}.
\end{align*}
Note that $\langle \varphi_{\mu}, \varepsilon_{n} \rangle$ has value 1 for $\mu=(-1,\cdots, -1)$ and zero for the others. Thereby, we obtain
$\langle e_{\lambda},\varepsilon_{n}\rangle=\sum_{\mu \in \mathcal{DP}(n)} u_{\lambda \mu} \langle \varphi_{\mu}, \varepsilon_{n} \rangle=u_{\lambda,(-1,\cdots, -1)}$. Eventually, we have $u_{\lambda,(-1,\cdots, -1)}=(-1)^{|S_{\hat{\lambda}}|}\cdot \frac{|\mathcal{C}(\lambda)|}{|W_{n}|}$.
\end{Prf}

Notice that calculation of the inner product $\langle e_{\lambda}, 1_{W_{n}}\rangle$ leads to the following corollary.

\begin{Cor}
Let $\lambda \in \mathcal{DP}(n)$. Then
\[
|W_{n}|\sum_{\mu \in \mathcal{DP}(n)}u_{\lambda,\mu}=|\mathcal{C}(\lambda)|.
\]
\end{Cor}

By means of this Corollary and the matrices $(u_{\lambda \mu})_{\lambda, \mu \in \mathcal{DP}(n)}$, we can readily determine the size of the conjugacy classes of $W_{n}$ . The proof of the next lemma is an immediate consequence of the inner product of characters.

For a subset $X$ of $W_{n}$, we denote by the subspace $\textrm{Fix}(X)=\{v \in \mathbb{R}^{n} : \forall x \in X~ x(v)=v \}$ of $\mathcal{R}^{n}$ fixed by $X$ and let write $W_{\textrm{Fix}(X)}=\{w \in W_{n} : \forall v \in \textrm{Fix}(X)~w(v)=v \}$ for the stabilizer of $\textrm{Fix}(X)$ in $W_{n}$. By [7], the set $W_{\textrm{Fix}(X)}$ is called the \textit{parabolic closure} of $X$ and it is denoted by $A(X)$. For a $w \in W_{n}$, if we take $X=\{w\}$ then we write $\textrm{Fix}(w)$ and $A(w)$ instead of $\textrm{Fix}(\{w\})$ and $A(\{w\})$, respectively. By [1], if $A(w)$ is $W_{n}$-conjugate to $W_{J}$ for some $J \subset S_{n}$, then we say that $w$ is of \textit{type} $J$.

\begin{Prop}\label{Proposition}
If $B \in \mathcal{SC}^{+}(n)$, then we have $A(W_{B})=W_{n}$.
\end{Prop}

\begin{Prf}
Since $B \in \mathcal{SC}^{+}(n)$, we have $\mathcal{T}_{n} \leq W_{B}$ and so $w_{n} \in W_{B}$. By considering $w_{n}$ as a linear map $-id_{\mathbb{R}^{n}}$, we obtain $Fix(w_{n})=\{\vec{0}\}$. Thus, the parabolic closure of $w_{n}$ is $A(w_{n})=W_{Fix(w_{n})}=W_{n}$. Because of the relation $w_{n} \in W_{B} \subset A(c_{B})=A(W_{B})$, we get $w_{n} \in A(c_{B})$. By [12], the inclusion $A(w_{n})\subset A(c_{B})=A(W_{B})$ holds. If we take into account the fact that $A(w_{n})=W_{n}$, then we have $A(W_{C})=W_{n}$. This completes the proof.
\end{Prf}

As a consequence of Proposition 3.6, if $B \in \mathcal{SC}^{+}(n)$, then the parabolic closure of $W_{B}$ is $W_{n}$ and each element of $\mathcal{C}(\boldsymbol{\lambda}(B))$ is of type $S_{n}$.

\begin{Lem}
Let $A$ be a signed composition of $n$. Then $w_{n}$ belongs to $W_{A}$ if and only if $A \in \mathcal{SC}^{+}(n)$.
\end{Lem}

\begin{Prf}
When $A$ is a positive signed composition of $n$, we can easily see from the proof of Proposition 3.6 that $w_{n}$ is an element of $W_{A}$. Conversely, let $w_{n}$ be in $W_{A}$. We suppose that $A=(a_{1}, \cdots, a_{i}, \cdots, a_{r})$ is not a positive signed composition of $n$. Then there exists $a_{i}<0$ for some $i,~1 \leq i \leq r$. Thus from the definition of $W_{A}$, we obtain $t_{|a_{1}|+\cdots+|a_{i}|} \not \in S_{A}^{'}$. Hence for $x \in W_{A}$ and $e_{|a_{1}|+\cdots+|a_{i}|} \in \mathbb{R}^{n}$, we have $x(e_{|a_{1}|+\cdots+|a_{i}|})=e_{|a_{1}|+\cdots+|a_{i}|}$ and so $e_{|a_{1}|+\cdots+|a_{i}|} \in Fix(W_{A})$. This is a contradiction, because  from Proposition 3.6 the subspace $\textrm{Fix}(W_{A})$ consists only $\vec{0}$. Therefore, we have $A \in \mathcal{SC}^{+}(n)$.
\end{Prf}

\begin{Prop}
If the set $\mathcal{C}(S_{n})$ denotes the set of all elements of type $S_{n}$, then we have
\begin{equation}
    \mathcal{C}(S_{n})=\coprod _{A \in \mathcal{SC}^{+}(n)}\mathcal{C}(\boldsymbol{\lambda}(A)).
\end{equation}
\end{Prop}

\begin{Prf}
From Proposition 3.6, for all $A \in \mathcal{SC}^{+}(n)$ every element of $\mathcal{C}(\boldsymbol{\lambda}(A))$ is of type $S_{n}$ and so the reverse inclusion holds. Now let $w \in \mathcal{C}(S_{n})$. Then $w$ is $W_{n}$-conjugate to $c_{A}$ for some $A\in \mathcal{SC}(n)$. Thus we get $A(w)=A(c_{A})=A(W_{A})=W_{n}$. From here, for every $x \in W_{n}$ and every $v \in Fix(W_{A})$ we obtain $x(v)=v$. In particular, if we take $w_{n}=-id_{\mathcal{R}^{n}} \in W_{n}$, then it is seen that $Fix(W_{A})$ includes just $\{ \vec{0} \}$. Thus $w_{n}$ is an element of  $\in W_{A}$. Otherwise, if $A \not \in \mathcal{SC}^{+}(n)$, then from the proof of Lemma 3.3 we get $Fix(W_{A}) \neq \{ \vec{0} \}$, which is a contradiction. Hence $A \in \mathcal{SC}^{+}(n)$. Thus $w$ belongs to $\mathcal{C}(\boldsymbol{\lambda}(A))$ and so it is seen that the inclusion $\mathcal{C}(S_{n})\subset \coprod _{A \in \mathcal{SC}^{+}(n)}\mathcal{C}(\boldsymbol{\lambda}(A))$ satisfies. It is required.
\end{Prf}

Since the exponents of $W_{n}$ are in turn $1, 3, \cdots, 2n-1$, then from \cite{Br1} the number of elements of $S_{n}$-type is equal to the product of exponents of $W_{n}$ and so $|\mathcal{C}(S_{n})|= 1\cdot 3 \cdots 2n-1$. By the equality (3.3), we obtain the formula
\begin{equation}
    |\mathcal{C}(S_{n})|=\sum_{A \in \mathcal{SC}^{+}(n)}|\mathcal{C}(\boldsymbol{\lambda}(A))|.
\end{equation}

Thus Proposition 3.7 gives us an alternative method to compute the number of elements of type $S_{n}$.

\begin{Lem}
Let $\varepsilon_{n}$ be the sign character of $W_{n}$ and let $A \in \mathcal{SC}(n)$. Then $\textrm{res}_{W_{A}}^{W_{n}}\varepsilon_{n}$ is an irreducible character of $W_{A}$.
\end{Lem}

Let $A \in \mathcal{SC}(n)$. When $A$ is not a parabolic signed composition of $n$, then $l_{A}(w)$ may not equal to $l(w)$ for some $w\in W_{A}$. The following lemma gives an interesting relation between these length functions.

\begin{Lem}
Let $A \in \mathcal{SC}(n)$. Let $l$ be the length function of $W_{n}$ attached to $S_{n}$ and let $l_{A}$ be the length function of $W_{A}$ in relation to $S_{A}$. Then for every $w \in W_{A}$
$$l(w)\equiv l_{A}(w)~(\textrm{mod}~2).$$
\end{Lem}

\begin{Prf}
We know that $t_{1}:=t$ and $t_{i}:=s_{i-1}t_{i-1}s_{i-1}$ for each i, $2 \leq i \leq n$. It is immediately seen that $l(t_{i})=2i-1$ according to the generating set $S_{n}$. Because $s_{i-1}s_{i-2} \cdots t_{1} \cdots s_{i-2}s_{i-1}$ is a reduced expression for $t_{i}$, then $l(t_{i})$ is equal to $2i-1$. Therefore for every i, $1 \leq i \leq n$, $l(t_{i})\equiv 1 (\textrm{mod}~2)$. Let $w\in W_{A}$. By virtue of $S_{A}\subset S'_{n}$, a reduced expression of $w$ in terms of the elements of $S_{A}$ can also include some elements of $T_{n}$ besides those of $S_{n}$. The lenght of $t_{i}$ which may be included in this reduced expression is 1 under $l_{A}$, and also $2i-1$ in terms of $l$. Thus, the length of $t_{i}$ is odd number with respect to both length functions. However, since the length of $s_{i}$ according to both $l_{C}$ and $l$ is $1$, then the statement of being odd or even for $w$ is the same in terms of both length functions. If we denote by $N_{A}(w)$ and $N(w)$ the sets of positive roots in $\Psi_{A}^{+}$ and $\Psi^{+}$ transformed into negative roots by $w$, respectively, then it follows that
$$N_{A}(w) \subset N(w).$$
It is seen from this fact that the reduced expression of $w\in W_{A}$ in terms of the fundamental reflections in $S_{A}$ is not longer than relative to those of $S_{n}$. Let $A=(a_{1}, \cdots, a_{r})$. Moreover, then because of $W_{A}\cong W_{a_{1}} \times \cdots \times W_{a_{r}}$, expanding the reduced expression of $w\in W_{A}$ according to simple reflections of $S_{A}$ by using the relations $t_{i}=s_{i-1}s_{i-2} \cdots t_{1} \cdots s_{i-2}s_{i-1}$ yield that reduced expression of $w$ in terms of simple reflections in $S_{n}$ can be obtained without any simplification. Hence $l(w)\equiv l_{C}(w)~(\textrm{mod}~2)$, as required.
\end{Prf}

As a result of the previous lemma, we get
$$\varepsilon_{n}(w)=(-1)^{l(w)}=(-1)^{l_{A}(w)}=\varepsilon_{A}(w).$$
Thus from Lemma 4 and 5, for every  $A \in \mathcal{SC}(n)$ it is obtained that $\textrm{res}_{W_{A}}^{W_{n}}\varepsilon_{n}=\varepsilon_{A}$.

\begin{Exam} \normalfont
For a concrete example, let $A=(-2, 3, -1, -3, 1) \in \mathcal{SC}(10)$, $S_{A}=\{s_{1}\} \cup \{t_{3}, s_{3}, s_{4}\} \cup \{s_{7}, s_{8}\} \cup \{t_{10}\} \subset S_{10}^{'}$ and $S_{A}=\{s_{1}\} \cup \{t_{3}, s_{3}, s_{4}, t_{4}, t_{5}\} \cup \{s_{7}, s_{8}\} \cup \{t_{10}\}$. Then $W_{A}\cong W_{\bar{2}} \times W_{3} \times W_{\bar{1}} \times W_{\bar{3}} \times W_{1}$. For $w=s_{7}t_{3}s_{3}s_{1}t_{10} \in W_{A}$, $l_{A}(w)=5$ and also
$$w=s_{7}t_{3}s_{3}s_{1}t_{10}=s_{7}s_{2}s_{1}t_{1}s_{1}s_{2}s_{3}s_{1}s_{9}s_{8}s_{7}s_{6}s_{5}s_{4}s_{3}s_{2}s_{1}t_{1}s_{1}s_{2}s_{3}s_{4}s_{5}s_{6}s_{7}s_{8}s_{9}\in W_{n},$$ so $l(w)=27$. It follows that $l(w)\equiv l_{A}(w) \equiv 1 (mod~2)$.
\end{Exam}

\begin{Theo}
Let $A \in \mathcal{SC}(n)$ and $\lambda \in \mathcal{DP}(n)$. Then
\[
\sum_{\mu \in \mathcal{DP}(n)}u_{\lambda \mu}a_{\hat{\mu} A (-1,\cdots,-1)}=(-1)^{|S_{\hat{\lambda}}|}\frac{|\mathcal{C}(\lambda)\cap W_{A}|}{|W_{A}|},
\]
where $a_{\mu A (-1,\cdots,-1)}=|\{x\in D_{\hat{\mu} A} : ~^{x^{-1}}\hat{\mu} \cap A=(-1,\cdots,-1)\}|$.
\end{Theo}

\begin{Prf}
The term $d_{(-1,\cdots,-1)}$ in the multiplication $d_{\hat{\mu}}d_{A}$ lies in the summand $\sum _{x\in D_{\hat{\mu}A}}d_{f_{\hat{\mu}A}(x)} $ from the structure of $\textrm{Ker}\theta_{n}$ and Proposition 2.1 (a). If we write the coefficient of $d_{(-1,\cdots,-1)}$ in this summand as $a_{\hat{\mu} A (-1,\cdots,-1)}$, and so we get
$$a_{\hat{\mu} A (-1,\cdots,-1)}=|\{ x \in D_{\hat{\mu}A} : f_{\hat{\mu}A}(x)=(-1,\cdots,-1)\}|.$$
By using Proposition 2.1 (a) along with the fact $f_{\hat{\mu}A}(x) \equiv_{A} {^{x^{-1}}\hat{\mu} \cap A}$, it is seen that there is the equivalence $^{x^{-1}}\hat{\mu} \cap A \equiv_{A} (-1,\cdots,-1)$. Since no elements in $\mathcal{SC}(n)$ is congruent to $(-1,\cdots,-1)$ except for $(-1,\cdots,-1)$, it then follows that $^{x^{-1}}\hat{\mu} \cap A = (-1,\cdots,-1)$. Hence we have deduced the equality $a_{\hat{\mu}A(-1,\cdots,-1)}=|\{ x \in D_{\hat{\mu}A} :~^{x^{-1}}\hat{\mu} \cap A=(-1,\cdots,-1)\}|$ holds. Therefore, by Frobenius Reciprocity and Mackey Theorem, we have
\begin{align*}
\langle \xi_{\lambda}, \textrm{ind}_{W_{A}}^{W_{n}}\varepsilon_{A} \rangle &=\sum_{\mu \in \mathcal{DP}(n)}u_{\lambda \mu} \sum_{x\in D_{\hat{\mu}A}} \langle \textrm{ind}_{W_{^{x^{-1}}\hat{\mu} \cap A}}^{W_{A}}1_{^{x^{-1}}\hat{\mu} \cap A}, \varepsilon_{A}\rangle \\
&=\sum_{\mu \in \mathcal{DP}(n)}u_{\lambda \mu} \sum_{\substack{ x\in D_{\hat{\mu}A} \\ ^{x^{-1}}\hat{\mu} \cap A=(-1,\cdots,-1)}} 1_{^{x^{-1}}\hat{\mu} \cap A} \\
&=\sum_{\mu \in \mathcal{DP}(n)}u_{\lambda \mu}a_{\hat{\mu} A (-1,\cdots,-1)}.
\end{align*}
Also, $\varepsilon(w)$ has the same value for every $w \in \mathcal{C}(\lambda)$ and so $\varepsilon(w)=\varepsilon(c_{\lambda})=(-1)^{|S_{\hat{\lambda}}|}$. Therefore, by Lemma 3.4, we have
\begin{align*}
\langle \xi_{\lambda}, \textrm{ind}_{W_{A}}^{W_{n}} \varepsilon_{A} \rangle &=\frac{1}{|W_{A}|} \sum_{w \in \mathcal{C}(\lambda) \cap W_{A}} (-1)^{l_{A}(w^{-1})}\\
&=\frac{1}{|W_{A}|} \sum_{w \in \mathcal{C}(\lambda) \cap W_{A}} (-1)^{l(w)}=\frac{1}{|W_{A}|}(-1)^{|S_{\hat{\lambda}}|}|\mathcal{C}(\lambda) \cap W_{A}|
\end{align*}
Putting these two results together, we see that theorem is proved.
\end{Prf}

\begin{Exam} \normalfont
For $n=2$, let $S_{2}=\{s,t\}$, $S_{2}^{'}=\{s,t, sts\}$ and
\begin{itemize}
\item [$\bullet$] $\mathcal{SC}(2)=\{(2),(1,1),(1,-1),(-,1),(-2),(-1,-1)\}$;
\item [$\bullet$] $\mathcal{DP}(2)=\{(2),(1,1),(1,-1),(-2),(-1,-1)\}.$
\end{itemize}
The conjugate classes of $W_{2}$ are as follows:
$\mathcal{C}((2))=\{st,ts\}$; $\mathcal{C}((1,1))=\{stst\}$; $\mathcal{C}((1,-1))=\{t,sts\}$;  $\mathcal{C}((-2))=\{s,tst\}$; $\mathcal{C}((-1,-1))=\{1\}.$
For $\lambda,\mu \in \mathcal{DP}(2)$
\begin{equation*}
(\varphi_{\lambda}(c_{\mu}))=
\begin{pmatrix}
  1 & 1 & 1 & 1 & 1 \\
  0 & 2 & 2 & 0 & 2 \\
  0 & 0 & 2 & 0 & 4 \\
  0 & 0 & 0 & 2 & 4 \\
  0 & 0 & 0 & 0 & 8 \\
\end{pmatrix},~
(u_{\lambda \mu})=
\begin{pmatrix}
  1 & \frac{-1}{2} & 0 & \frac{-1}{2} & \frac{1}{4} \\
  0 & \frac{1}{2} & \frac{-1}{2} & 0 & \frac{1}{8} \\
  0 & 0 & \frac{1}{2} & 0 & \frac{-1}{4} \\
  0 & 0 & 0 & \frac{1}{2} & \frac{-1}{4} \\
  0 & 0 & 0 & 0 & \frac{1}{8} \\
\end{pmatrix}.
\end{equation*}

\begin{itemize}
\item [$\bullet$] $|\mathcal{C}((2))|=|W_{2}|\cdot \sum_{\mu \in\mathcal{DP}(2)}u_{(2),\mu}
=8.(1+ \frac{-1}{2} + 0 + \frac{-1}{2} + \frac{1}{4})=2$;
\end{itemize}
For $C=(2),(1,1),(1,\bar{1}),(\bar{2}),(\bar{1},\bar{1})\in \mathcal{DP}(2)$, the Coxeter generating sets are $S_{C}=S_{2}=\{s,t\}$,~$\{t,stst\}$,~$\{t\}$,~$\{s\}$,~$\emptyset$ respectively. Then we have
\begin{align*}
u_{(2),(-1,-1)} &=\frac{(-1)^{|S_{2}|}}{2 \cdot 2}=\frac{1}{4} \\
u_{(1,1),(-1,-1)} &=(-1)^{|S_{(1,1)}|}\frac{|\mathcal{C}((1,1))|}{|W_{2}|}=\frac{1}{8} \\
u_{(1,-1),(-1,-1)} &=(-1)^{|S_{(1,-1)}|}\frac{|\mathcal{C}((1,-1))|}{|W_{2}|}=\frac{-1}{4} \\
u_{(-2),(-1,-1)} &=(-1)^{|S_{(-2)}|}\frac{|\mathcal{C}((-2))|}{|W_{2}|}=\frac{-1}{4} \\
u_{(-1,-1),(-1,-1)} &=(-1)^{|S_{(-1,-1)}|}\frac{|\mathcal{C}((-1,-1))|}{|W_{2}|}=\frac{1}{8}.
\end{align*}
Note that these results are the same as  from top to bottom entries in the last column of the matrix $(u_{\lambda, \mu})$.
As considering Corollary 3.1 and Proposition 3.7, we obtain that the set of all $S_{2}$-type elements is $\mathcal{C}(S_{2})=\mathcal{C}((2))\uplus\mathcal{C}((1,1))$ and so from (3.4) the number of type $S_{2}$ is $|\mathcal{C}(S_{2})|=|\mathcal{C}((2))|+|\mathcal{C}((1,1))|=3$. Also note that here $|\mathcal{C}(S_{2})|$ is equal to the product of exponents of $W_{2}$.

\end{Exam}

\begin{center}
ACKNOWLEDGEMENTS
\end{center}

The authors are indebted to Götz Pfeiffer for the opportunity to work together and for many useful discussions.

\section*{\normalsize{4- References}}


\vspace*{-1cm}

\begin{thebibliography}{99}


\bibitem{Br1} \textbf{  Bergeron, F., Bergeron, N., Howlett, R.B., Taylor, D.E.}, {\em A Decomposition of the Descent Algebra of a Finite Coxeter Group}, Journal of Algebraic Combinatorics,  1(1), 23-44, 1992.
\bibitem{Br2} \textbf{  Bonnafé, C., Hohlweg, C.}, {\em Generalized descent algebra and construction of irreducible characters of hyperoctahedral groups}, Ann. Inst. Fourier (Grenoble) 56 (2006), no. 1, 131-181.
\bibitem{Br3} \textbf { Bonnafé, C.}, {\em Representation theory of Mantaci-Reutenauer algebras}, Algebras and Representation Theory 11 (2008), no. 4, 307-346.
\bibitem{Br4} \textbf { Carter, R.W.}, {\em Conjugacy Classes in the Weyl Groups}, Compositio Math., 25 (1972), 1-59.
\bibitem{Br5} \textbf {  Curtis, C.W., Reiner, I.}, {\em Representation Theory of Finite Groups and Associative Algebras}, John Wiley and Sons, New York, 1962.
\bibitem{Br6} \textbf {Curtis, C.W., Reiner, I.}, {\em Methods of Representation Theory with Applications to Finite Groups and Orders}, Vol. II, John Wiley and Sons, 1987.
\bibitem{Br7} \textbf{ Douglass, J.M.; Pfeiffer, G.; Röhrle, G.} {\em On reflection subgroups of finite Coxeter groups}, Comm. Algebra 41 (2013), no. 7, 2574-2592.
\bibitem{Br8} \textbf{ Fleischmann, P.}, {\em On pointwise conjugacy of distinguished coset representatives in Coxeter groups}, J. Group Theory, 5 (2002), 269-283.
\bibitem{Br9} \textbf{ Geck, M., Pfeiffer, G.}, {\em Characters of Finite Coxeter Groups and Iwahori-Hecke Algebras}, London Mathematical Society Monographs, New Series, vol. 21, The Clarendon Press, Oxford University Press, New York, 2000.
\bibitem{Br10} \textbf{ Humphreys, J.E.}, {\em Reflection Groups and Coxeter Groups}, Cambridge Studies in  Advenced Math., vol. 29,  Cambridge University Press, 1990.
\bibitem{Br11} \textbf{ Mantaci, R., Reutenauer, C. }, {\em A generalization of Solomon's algebra for hyperoctahedral groups and other wreath products}, Comm. Algebra 23 (1995), no. 1, 27-56.
\bibitem{Br12}\textbf { Solomon, L.}, {\em A Mackey formula in the group ring of a Coxeter group}, J. Algebra 41 (1976), no. 2, 255-264.





\end{thebibliography}
\end{document}